\documentclass{article}
\usepackage{amsmath}
\usepackage{amsthm}
\usepackage{latexsym}
\usepackage{amssymb}
\usepackage[all]{xy}
\usepackage{graphicx}

\theoremstyle{definition}
\newtheorem{thm}{Theorem}[section]
\newtheorem{defi}[thm]{Definition}
\newtheorem{prop}[thm]{Proposition}
\newtheorem{cor}[thm]{Corollary}

\newtheorem{lemma}[thm]{Lemma}
\newtheorem{rem}[thm]{Remark}

\newcommand\PP{\ensuremath{\mathbb{P}}}
\newcommand\ZZ{\ensuremath{\mathbb{Z}}}
\newcommand\CC{\ensuremath{\mathbb{C}}}

\newcommand\RR{\ensuremath{\mathbb{R}}}

\newcommand\vect[1]{\mbox{\boldmath$#1$}}

\title{The totally nonnegative part of the finite Toda lattice via a reducible rational curve}
\author{
Shinsuke Iwao${}^1$\footnote{E-mail: iwao@gem.aoyama.ac.jp}, 
Kyo Nishiyama${}^2$\footnote{E-mail: kyo@gem.aoyama.ac.jp}
and 
Noboru Ogawa${}^3$\footnote{E-mail: nogawa@tsc.u-tokai.ac.jp}   \\[5pt]
${}^{1,2}$College of Science and Engineering, Aoyamagakuin University,\\
5-10-1 Fuchinobe, Cyu\={o}-ku, Sagamihara-shi, Kanagawa Japan 252-5258\\
${}^3$Department of Mathematics, Tokai University,\\
4-1-1 Kitakaname,
Hiratsuka-shi, Kanagawa Japan 259-1292
}

\begin{document}
\maketitle

{\abstract
A totally nonnegative matrix is a real-valued matrix whose minors are all nonnegative. 
In this paper, we concern with the totally nonnegative structure of the finite Toda lattice, a classical integrable system, which is expressed as a differential equation of square matrices.
The Toda flow naturally translates into a (multiplicative) linear flow on the (generalized) Jacobi variety associated with some reducible rational curve $X$.
This correspondence provides an algebro-geometric characterization of the totally positive part of the Toda equation.
We prove that the totally nonnegative part of the finite Toda lattice is isomorphic to a connected component of $\mathrm{Jac}(X)_\RR$, the real part of the generalized Jacobi variety $\mathrm{Jac}(X)$, as semi-algebraic varieties.
}

\section{Introduction}\label{sec1}

A real matrix is said to be {\it totally positive} ({\it resp.\,totally nonnegative}) if all its minors are positive ({\it resp.\,nonnegative}).
In this article, we study the {\it totally nonnegative part} of the phase space of the {\it finite Toda lattice}, that is one of the most basic classical integrable system.
Since the Toda lattice was originally discovered as a physical system, it would be natural to expect that the totally nonnegative part possesses all the essential structures of the system, although it is ``mathematically natural'' to consider the phase space of complex numbers.
We give an algebro-geometric characterization of the totally nonnegative part of the finite Toda lattice via the Krichever construction 
associated with a singular algebraic curve.

\subsubsection*{The aim}

The methods to construct solutions of the finite Toda lattice were mainly developed in 1970's by many authors. 
It is well-known that there exist quite a few approaches to give explicit solutions such as the {\it inverse scattering method} \cite{PhysRevB.9.1924}, the method of the {\it QR decomposition} \cite{Kostant1979195,Olshanetsky1979}, the method of {\it Wronskian matrix solutions}, the {\it continued fraction expansion method} \cite{Moser1975}, {\it etc.}

In this study, we use so-called {\it Krichever construction} associated with some singular curve \cite{krichever1977methods,Krichever2002}. 
In 2002, Krichever and Vaninsky \cite{Krichever2002,Vaninsky2003283} pointed out the relation between the finite Toda lattice and some kind of reducible curve, and found the fact that the Baker-Akhiezer function associated with that curve gives rational solutions of the finite Toda lattice.
After that, Sklyanin \cite{sklyanin2013bispectrality} obtained more explicit expressions of these rational solutions in terms of the {\it dual variable} by making more detailed analysis on the reducible curve.

Inspired by their works \cite{Krichever2002,sklyanin2013bispectrality}, we give an geometric interpretation of the totally non-negative part of the finite Toda lattice by re-interpreting their construction as a {\it generalized Jacobi inversion problem}.
We are also inspired by the paper \cite{kodama2014} in 2014 by Kodama and Williams, who have investigated the combinatorial structure of the totally nonnegative part of the system.
(In \cite{kodama2014}, they have also studied the full Kostant-Toda hierarchy, which is a generalization of the Toda lattice.
It should be an interesting problem to give a similar characterization of the totally nonnegative part of it.
Also see \cite{kodama2008finite} for a detailed survey of the topological properties of the finite Toda lattice.)

The outline of the result of this paper is as follows:
Let $L\in M_N(\CC)$ be a complex matrix of size $N$.
Its characteristic polynomial $\det(\lambda E-L)=\lambda^N-H_1(L)\lambda^{N-1}+H_2(L)\lambda^{N-2}-\cdots+(-1)^NH_N(L)$ defines the algebraic map
\[
\begin{array}{cccccccccc}
 M_N(\CC) &\xrightarrow{H}  & \CC^N & \xleftarrow{\gamma}& \CC^N/\mathfrak{S}_N\\[10pt]
L & \mapsto  & (H_1(L),\dots,H_N(L)) \\
  & & (e_1(\Lambda),\dots,e_N(\Lambda)) & \leftarrow& \Lambda=\{\lambda_1,\dots,\lambda_N\},
\end{array}
\]
where $e_i(\Lambda)$ is the $i$-th symmetric function in $\lambda_1,\dots,\lambda_N$.
Let $\Gamma$ (Eq.\,(\ref{eq0})) be the phase space of the finite Toda lattice (Eq.\,(\ref{eq1})). 
The subset $\mathcal{T}_\Lambda:=H^{-1}(\gamma(\Lambda))\cap \Gamma$ is called an {\it isospectral set}.
The Toda equation determines a time-independent flow (so-called the {\it Toda flow}) on each isospectral set $\mathcal{T}_\Lambda$.
For any $\Lambda$, the isospectral set $\mathcal{T}_\Lambda$ is connected and of complex dimension $N-1$ as a subvariety of $M_N(\CC)$.

We consider the set $\mathcal{T}_{\Lambda}^\geq:=\mathcal{T}_\Lambda\cap \{\mbox{totally non-negative matrices}\}$, the {\it totally non-negative part of} $\mathcal{T}_\Lambda$.
It is known that, if $\mathcal{T}_{\Lambda}^\geq$ is non-empty, the numbers $\lambda_1,\dots,\lambda_N\in \Lambda$ must be distinct, real and positive. 
(Theorem \ref{gmt}).

The following is the main theorem of this paper:
\begin{thm}[Theorem \ref{thm3.8}]
Assume $\mathcal{T}_\Lambda^{\geq}\neq \emptyset$.
Then, the image of $\mathcal{T}_\Lambda^\geq$ by the {\it linearization map} $\Phi:\mathcal{T}_\Lambda\to \mathrm{Jac}(X)$ (Eq.\,(\ref{kric})) is described as follows:
\[
\Phi(\mathcal{T}_\Lambda^{\geq})=
\left\{
[F_1:F_2:\dots:F_N]\in\mathrm{Jac}(X)\,\left\vert\, (-1)^iF_i>0\mbox{ for all } i
\right.
\right\},
\]
where $X$ is a reducible curve defined in \S \ref{sec3}, and $\mathrm{Jac}(X)\simeq (\CC^\times)^N/\CC^\times$ is the (generalized) Jacobi variety associated with $X$.
%
\end{thm}

Especially, the image $\Phi(\mathcal{T}_\Lambda^{\geq})$ is a positive cone of (real) dimension $N-1$.
The Abel-Jacobi map $\Phi$ induces the isomorphism
\[
\Phi\vert_{\mathcal{T}_\Lambda^{\geq}}: \mathcal{T}_\Lambda^\geq\to \mathrm{Jac}(X)_\RR^0
\]
of semi-algebraic varieties, where $\mathrm{Jac}(X)_\RR^0$ is the connected component of $\mathrm{Jac}(X)_\RR$ (\S \ref{sec3.3}) which contains the identity element $[1:1:\dots:1]$.

This result provides several properties of the Toda flow.
For example, the totally non-negative part $\mathcal{T}_\Lambda^\geq$ is connected and closed under the Toda flow.
Any orbit of a totally non-negative matrix never contain a blowup point.

\subsubsection*{Contents of the paper}

The contents of the paper is as follows:

In Section \ref{sec2}, we demonstrate the algebro-geometric construction of determinant solutions of the Toda lattice by solving the (generalized) Jacobi inversion problem associated with a singular spectral curve.
The singular curve which we use here is the same thing as in the papers \cite{Krichever2002,Vaninsky2003283}. 
Although the conclusion (Theorem \ref{kai1}) of this section is not new at all, our method provides more suitable geometric interpretations of some datum related to the curve.
(For example, the ``gluing condition'' in \cite[\S 3]{Krichever2002} can be understood as an element of $\mathrm{Jac}(X)$ in our context).

In Section \ref{sec3a}, we give a proof of the main theorem (Theorem \ref{thm3.8}).
Our proof (\S \ref{sec3.2}) essentially depends on the explicit form of the determinant solution.

\subsubsection*{Remark: two meanings of ``totally non-negativity''.}

When referring to the ``totally non-negativity'' of the Lax equation $\dfrac{d}{dt}L=[L,L_-]$, we should be careful not to confuse the following two meanings:
\begin{enumerate}
\item the totally non-negativity of the {\it companion matrix} of  $L$.
\item the totally non-negativity of $L$ itself.
\end{enumerate}
Many researchers maybe familiar with the first one, that has been well-studied in researches on blowup solutions of the Toda equation. 
(The ``TNN part'' of Kodama-Williams belongs to this.)
On the other hand, among researches on discrete integrable systems and their tropicalizations, the totally non-negativity in the second meaning is a main object to be studied.
It is straightforward to check that the totally non-negativity in the second meaning implies the first one.

In this study, we always refer to the totally non-negativity in the second meaning.

\section{The solution of the finite Toda lattice in terms of singular theta functions}\label{sec2}

\subsection{The finite Toda lattice in Lax formalism}

Let $N$ be a positive integer.
Put $G=GL_N(\CC)$ and $\mathfrak{g}=\mathfrak{gl}_N(\CC)$.
Define the subgroups $B^+,U^-\in G$ by
\begin{gather*}
B^+:=\{x\in G\,\vert\,\mbox{upper triangular}\},\qquad
U^-:=\{x\in G\,\vert\,\mbox{lower triangular with diagonals $1$}\}
\end{gather*}
and the subalgebras $\mathfrak{g}_\pm\subset\mathfrak{g}$ by
\begin{gather*}
\mathfrak{g}_+:=\{x\in \mathfrak{g}\,\vert\,\mbox{upper triangular}\}=\mathrm{Lie}(B^+),\qquad
\mathfrak{g}_-:=\{x\in \mathfrak{g}\,\vert\,\mbox{strictly lower triangular}\}=\mathrm{Lie}(U^-).
\end{gather*} 
It is well-known that the product map $B^+\times U^-\to G$; $(b,u)\mapsto b\cdot u$ is an open embedding so that we have $\mathfrak{g}=\mathfrak{g}^+\oplus \mathfrak{g}^-$.
For $X\in\mathfrak{g}$, denote by $X=X_++X_-$ ($X_+\in \mathfrak{g}^+$, $X_-\in \mathfrak{g}^-$) the decomposition of $X$ along the direct sum decomposition.

Let $\Gamma\subset \mathfrak{g}$ be a subset defined by
\begin{equation}\label{eq0}
\Gamma:=
\left\{
\left(
\begin{array}{cccc}
a_1 & 1 &  &   \\ 
b_1 & a_2 & \ddots &    \\ 
 & \ddots  & \ddots & 1  \\ 
 &  &  b_{N-1} & a_N
\end{array} 
\right)
\in \mathfrak{g} \ ;\ 
b_1b_2\cdots b_{N-1}\neq 0
\right\}.
\end{equation}
In this section, we consider the {\it finite Toda lattice} \cite{toda1967vibration,toda1967wave}:
\begin{equation}\label{eq1}
\frac{d}{dt}L=[L, L_-],\qquad L=L(t)\in \Gamma.
\end{equation}
Fix an initial value $L_0=L(0)\in\Gamma$ of the system (\ref{eq1}).
For a real number $t$, $\exp(tL_0)=a(t)^{-1}b(t)\in B^+\cdot U^-$ be the LU decomposition.

\begin{lemma}[\cite{Symes1982275}]
$L(t):=a(t)L_0a(t)^{-1}=b(t)L_0b(t)^{-1}$ solves the equation (\ref{eq1}). $\qed$
\end{lemma}

\subsubsection{The Toda flow on the spectrum set}\label{sec2.1.1}

From the Lax form (\ref{eq1}), the characteristic polynomial 
\begin{equation}\label{char}
f(\lambda):=(-1)^N\det(L(t)-\lambda\cdot E)=(-1)^N\det(L_0-\lambda\cdot E)
\end{equation}
is independent of $t$.
Since we are interested in the totally non-negative structure, we restrict ourselves to the case when all the roots of $f(\lambda)$ are simple (see Theorem \ref{gmt}).
We denote these roots by $\lambda_1,\dots,\lambda_N\in \CC$.
Let $C$ be the finite set $C:=\{\lambda_1,\dots,\lambda_N\}\subset \CC$, the spectrum set of $L(t)$.
The ring $\mathcal{O}:=\CC[\lambda]/(f(\lambda))$ is regarded as the ring of algebraic functions over $C$.
By the Chinese remainder theorem, $\mathcal{O}$ admits the ring isomorphism
\begin{equation}\label{eq3}
\mathcal{O}\xrightarrow{\sim} \CC^N;\qquad \varphi\mapsto (\varphi(\lambda_1),\dots,\varphi(\lambda_N)).
\end{equation}

For each $t\in \RR$, consider a column vector $\vect{v}(t):=(v_1(t),\dots,v_N(t))^T$, ($v_i(t)\in\mathcal{O}$) which satisfies the linear equation
\begin{equation}\label{eq4}
(L(t)-\lambda E)\vect{v}(t)\equiv \vect{0}\mod{f(\lambda)}.
\end{equation}

The cofactor expansion of $(L(t)-\lambda E)$ along the $N$-th row yields the column vector
\[
\vect{v}_-(t)=((-1)^{N+j}\Delta_{N,j}(t;\lambda))_{j=1}^N,\qquad \mbox{where }\Delta_{i,j} \mbox{ is a $(i,j)$-minor of $(L(t)-\lambda E)$}.
\]
As a polynomial in $\lambda$, the $i$-th component of $\vect{v}_-(t)$ is expressed as $\lambda^{i-1}+O(\lambda^{i-2})$.
On the other hand, the expansion along the first row yields the column vector 
$$
\vect{v}_+(t)=((-1)^{1+j}\Delta_{1,j}(t;\lambda))_{j=1}^N
$$
whose $i$-th components is of the form $(b_1b_2\cdots b_{i-1})\lambda^{N-i}+O(\lambda^{N-i-1})$.
Since the matrix $(L(t)-\lambda E)$ is always expressed as $P\cdot \mathrm{diag}(1,\dots,1,f(\lambda))\cdot Q$, ($P,Q$ are invertible matrix over $\CC[\lambda]$), the two vectors $\vect{v}_-(t)$ and $\vect{v}_+(t)$ are in proportion:
\begin{equation}\label{eq7}
F(t)\cdot \vect{v}_{-}(t)= \vect{v}_{+}(t),\qquad \exists ! F(t)\in \mathcal{O}^\times.
\end{equation}
Comparing the first components on the both sides, we have
$
F(t)=\Delta_{1,1}(t;\lambda).
$
Therefore, $F(t)$ is nothing but the {\it chop integral} introduced in \cite{kodama2014,kodama2008finite}.

The time dependence of $F(t)$ is calculated as follows.
From $L(t)=a(t)L_0a(t)^{-1}=b(t)L_0b(t)^{-1}$, we have two equations
\begin{gather*}
L(t)\{a(t)\vect{v}(0)\}=\lambda \{a(t)\vect{v}(0)\}, \quad L(t)\{b(t)\vect{v}(0)\}=\lambda \{b(t)\vect{v}(0)\}. 
\end{gather*}
Therefore, all the vectors $\vect{v}_{\pm}(t)$, $a(t)\vect{v}_{\pm}(0)$, $b(t)\vect{v}_{\pm}(0)$ are in proportion.
From the facts $a(t)\in U^-$ and $b(t)\in B^+$, we conclude
\begin{enumerate}
\item $\vect{v}_-(t)=a(t)\vect{v}_-(0)$ and 
\item $\vect{v}_+(t)=k(t)b(t)\vect{v}_+(0)$ for some non-zero complex number $k(t)\in \CC^\times$.
\end{enumerate}
Combining these equations and (\ref{eq7}), we obtain $F(t)\cdot a(t)\vect{v}_-(0)=k(t)b(t)\vect{v}_+(0)$.
Therefore, $F(t)\cdot\vect{v}_-(0)=k(t)a(t)^{-1}b(t)\vect{v}_+(0)
 =k(t) \exp(tL_0)\vect{v}_+(0)
 =k(t)\exp(t\lambda)\vect{v}_+(0)$.
This yields 
\begin{equation}
F(t)=k(t)\exp(t\lambda)\cdot F(0),\qquad F(t)\in \mathcal{O}^\times
\end{equation}
or in other words,
\begin{equation}\label{eq8}
\overline{F}(t)=\exp(t\lambda)\cdot \overline{F}(0),\qquad \overline{F}(t)\in \mathcal{O}^\times/\CC^\times,
\end{equation}
where $\overline{F}(t)=F(t)\,\mathrm{mod}\,{\CC^\times}$.




\subsection{The generalized Jacobi variety of the singular spectral curve.}\label{sec3}

Denote $\PP:=\PP^1(\CC)$. 
Let $X_-$ and $X_+$ be a pair of copies of $\PP$ with coordinate variables $x$ and $y$, respectively.
Define the reduced, reducible and nodal curve $X$ by gluing $X_-$ and $X_+$ along $\{x=\lambda_i\}\in X_-$ and $\{y=\lambda_i\}\in X_+$, ($i=1,2,\dots,n$) transversally.
Let $P_1,\dots,P_N\in X$ be the singular points (nodes) of $X$ and $\infty_+,\infty_-\in X$ be the points at infinity on $X_-$ and $X_+$, respectively\footnote{The affine part $X\setminus \{\infty_\pm\}$ is isomorphic to $\mathrm{Spec}\,[\CC[\lambda,\mu]/(\mu^2-\prod_i(\lambda-\lambda_i)^2)]$.}.
Later we will assign this singular curve $X$ with spectral datum of the finite Toda lattice.

\subsubsection{Definition of the generalized Abel-Jacobi map}

Let $H(X)$ be the group
$$
H(X)
:=\left\{(f,g)\,
\left\vert\,
\begin{array}{l}
\mbox{(i) $f$ ({\it resp.\,}$g$) is a rational function over $X_-$ ({\it resp.\,}$X_+$)},\\
\mbox{(ii) $f$ and $g$ have no poles and zeros at $P_1,\dots,P_N$},\\
\mbox{(iii) } f(P_n)=g(P_n),\quad
n=1,\dots,N
\end{array}
\right\}\right.,
$$
where the product is given by $(f,g)\cdot (f',g'):=(ff',gg')$.

Let $\mathrm{Div}(X):=\bigoplus_{p\in X\setminus\{P_1,\dots,P_N\}}{\ZZ\cdot p}$ be the divisor group of $X$ and
$\mathrm{Div}^d(X)$ be the subset of divisors of degree $d\in \ZZ$.
Further, define the subset $\mathrm{Div}^{a,b}(X)\subset \mathrm{Div}^{a+b}(X)$ by
\[
\mathrm{Div}^{a,b}(X):=\{D\in \mathrm{Div}^{a+b}(X)\,\vert\, 
\deg{(D\cap X_-)}=a,\ \deg{(D\cap X_+)}=b
\}.
\]
It follow that $\mathrm{Div}^d(X)=\coprod_{a+b=d}{\mathrm{Div}^{a,b}(X)}$.

For a pair of non-zero rational functions $v=(f,g)$, ($f$ ({\it resp.\,}$g$) is a rational function over $X_-$ ({\it resp.\,}$X_+$)), define the divisor $(v)\in \mathrm{Div}^{0,0}(X)$ by
$
(v):=(\mbox{the zeros of }v)-(\mbox{the poles of }v)=(v)_0-(v)_\infty
$.

For a divisor $D\in \mathrm{Div}^{0,0}(X)$, there exists a pair of rational functions $(f,g)$ whose divisor is $D$.
This determines the homomorphism 
\begin{equation}\label{eq9}
\textstyle
\mathrm{Div}^{0,0}(X)\to \mathcal{O}^\times/\CC^\times \simeq (\CC^\times)^N/\CC^\times;\qquad D\mapsto 
\left[\frac{g(\lambda_1)}{f(\lambda_1)}:\cdots:\frac{g(\lambda_N)}{f(\lambda_N)}
\right].
\end{equation}
Define the {\it Picard group} $\mathrm{Pic}^{a,b}(X)$ {\it of} $X$ of degree $(a,b)$ by
\[
\mathrm{Pic}^{a,b}(X):=\mathrm{Div}^{a,b}(X)/\sim,
\qquad D_1\sim D_2\iff D_1-D_2=(v),\ \ \exists v\in H(X).
\]
The homomorphism (\ref{eq9}) induces the map
$
J: \mathrm{Pic}^{0,0}(X)\to \mathcal{O}^\times/\CC^\times
$.

\begin{defi}
We define the {\it (generalized) Jacobi variety of} $X$ as $\mathrm{Jac}(X):=\mathcal{O}^\times/\CC^\times$.
The map $J:\mathrm{Pic}^{0,0}(X)\to \mathrm{Jac}(X)$ is called the {\it (generalized) Abel-Jacobi map}.
\end{defi}

\begin{prop}
$J$ is a group isomorphism.
\end{prop}
\proof 
Assume $J(D)=[1:\dots:1]$. 
Then, there exists some $v=(f,g)\in H(X)$ such that $D=(v)$, which implies $D\equiv 0$.
Therefore, $J$ is injective.
On the other hand, it is straightforward to prove $J$ to be surjective because there must exist a polynomial $h(\lambda)$ such that $[h(\lambda_1):\dots:h(\lambda_N)]=F$ for any $F=[F_1:\dots:F_N]\in \mathcal{O}^\times/\CC^\times$.
\qed

\subsubsection{General divisors}

Consider the subset
\[
S_{a,b}:=\{p_1+\dots+p_a+q_1+\dots+q_b\in \mathrm{Div}^{a,b}(X)
\,\vert\,p_i\in X_-,q_i\in X_+  \},\qquad (a,b\in\ZZ_{\geq 0})
\]
of formal sums of $(a+b)$ points.
We demonstrate the analog of the classical {\it Jacobi inversion problem}, that asks whether the map
$$
\phi_{a,b}:S_{a,b}\to \mathrm{Jac}(X);\qquad 
D\mapsto  J(D-a\cdot \infty_--b\cdot \infty_+)
$$
admits an inverse map.

\begin{prop}
If $a+b=N-1$, the map $\phi_{a,b}$ is surjective.
Moreover, for an element $A=[z_1:\dots:z_N]$ in a generic position of the Jacobi variety $\mathrm{Jac}(X)$, the inverse image $\phi_{a,b}^{-1}(A)$ consists of one element.
\end{prop}
\proof
The existence of the inverse image of $[z_1:\dots:z_N]\in\mathrm{Jac}(X)$ is equivalent to the existence of the pair of polynomials $f,g$ with $\dfrac{g(\lambda_n)}{f(\lambda_n)}=z_n$, $\deg f=a$ and $\deg g=b$.
Since the total number of coefficients of polynomials $f,g$ is $a+b+2=N+1$, we can conclude the claim by counting the number of variables. \qed

By this proposition, if $a+b=N-1$, we have the relation
\begin{equation}\label{eq11}
D\sim D',\ \ D'\in S_{a,b},\ \iff \phi_{a,b}(D)=\phi_{a,b}(D')\stackrel{\scriptsize\mbox{(if $D$ is in general position)}}{\iff } D=D'.
\end{equation}
In general, an element $D\in S_{a,b}$ with the property (\ref{eq11}) is called a {\it general divisor}.

\subsubsection{Theta functions}\label{sec5}

For distinct complex numbers $\lambda_1,\dots,\lambda_N\in \CC$, we define the multivalued functions $\theta_k(\,\cdot\,;\lambda_1,\dots,\lambda_N):(\CC^\times)^N\to \CC$ for $k=0,1,\dots,N$ by the following formula:
\begin{equation}\label{theta}
\theta_k(Z_1,\dots,Z_N;\lambda_1,\dots,\lambda_N):=
\frac{1}{\sqrt{Z_1\cdots Z_N}}
\det
\left(
\begin{array}{ccc@{\ }c@{\ }c|cc@{\ }c@{\ }c}
1 & \lambda_1 & \lambda_1^2 &  & \lambda_1^{N-k-1} & Z_1 & Z_1\lambda_1 &  & Z_1\lambda_1^{k-1}  \\ 
1 & \lambda_2 & \lambda_2^2 & \cdots & \lambda_2^{N-k-1} & Z_2 & Z_2\lambda_2 & \cdots &  Z_2\lambda_2^{k-1}\\ 
\vdots & \vdots & \vdots & \cdots & \vdots & \vdots &  \vdots & \cdots &  \vdots\\ 
1 & \lambda_N & \lambda_N^2 &  & \lambda_N^{N-k-1} & Z_N & Z_N\lambda_N  &  & Z_N\lambda_N^{k-1}
\end{array} 
\right).
\end{equation}
We call $\theta_k$ the {\it degenerated theta function}.

For the reducible curve $X$ defined in previous sections, we refer to the function $\theta_k(Z)=\theta_k(Z;\lambda_1,\dots,\lambda_k)$ as the {\it theta function associated with $X$}.
The following lemma provides the reason why it can be called the ``theta function''.
\begin{lemma}
For $[z_1:\dots:z_N]\in \mathrm{Jac}(X)$, we have 
$$
\theta_k(z_1,\dots,z_N)=0 \iff \exists D\in S_{k-1,N-k-1},\quad
\phi_{k-1,N-k-1}(D)=[z_1:\dots:z_N].
$$
The condition on the right hand side does not depend on the choice of the lift $(z_1,\dots,z_N)\in (\CC^\times)^N$.
$\qed$
\end{lemma}

For elements $\vect{A}=(A_1,\dots,A_N)$ and $\vect{B}=(B_1,\dots,B_N)$ of $(\CC^\times)^N$, define $\vect{A}\vect{B}=(A_1B_1,\dots,A_NB_N)$ and $\frac{\vect{A}}{\vect{B}}=(\frac{A_1}{B_1},\dots,\frac{A_N}{B_N})$.
Let $\varphi:X\setminus \{P_1,\dots,P_N\}\to (\CC^\times)^{N}$ be the continuous map
\[
\varphi(p):=
\left\lbrace 
\begin{array}{cc}
\left(\lambda_1-\lambda(p),\dots,\lambda_N-\lambda(p)
\right)^{-1},
 & p\in X_-\\
\left(\lambda_1-\lambda(p),\dots,\lambda_N-\lambda(p)
\right), & p\in X_+
\end{array}
\right..
\]
For a fixed element $ \vect{Z}\in (\CC^\times)^N$, define the multivalued function $\Theta_k(\ \cdot\ ;\vect{Z}):X\setminus\{P_1,\dots,P_N\}\to \CC$ by
$$
\Theta_k(p;\vect{Z}):= 
\left\{
\begin{array}{cc}
\theta_{k-1}(\vect{Z}/ \varphi(p)), & (p\in X_-)\\
\theta_{k}(\vect{Z}/ \varphi(p)), & (p\in X_+) 
\end{array}
\right..
$$
$\Theta_k(p;\vect{Z})$ ramifies around $p=P_i$ (the ramification number $=2$) and around $\infty_\pm$ (the ramification number $=2$ if $N$ is odd, and $=1$ if $N$ is even).

\begin{prop}[Analog of Riemann's vanishing theorem]
For $1\leq k\leq N$ and $\vect{Z}=(Z_1,\dots,Z_N)\in(\CC^\times)^N$, fix a divisor $D\in S_{k-1,N-k} $ with $\phi_{k-1,N-k}(D)=[Z_1:\dots:Z_N]$.
Then, the following (i--ii) hold:\\
(i) $\Theta_k(\,\cdot\,;\vect{Z})$ is identically $0$ if $D$ is non-general.\\
(ii) $\Theta_k(\,\cdot\,;\vect{Z})$ has zeros only at $D\in S_{k-1,N-k}$ if $D$ is general.
$\qed$
\end{prop}

\begin{cor}\label{cor2.7}
Suppose $D_1\in S_{k-1,N-k}$ and $D_2\in S_{l-1,N-l}$ to be general divisors.
Let $\vect{Z},\vect{Y}\in (\CC^\times)^N$ be the lifts of
$$
\phi_{k-1,N-k}(D_1),\ \phi_{l-1,N-l}(D_2)\in \mathrm{Jac}(X)\simeq (\CC^\times)^N/\CC^\times.
$$
Then, the ratio
$
\dfrac{\Theta_k(p;\vect{Z})}{\Theta_{l}(p;\vect{Y})}
$
is a single-valued function over $X\setminus \{P_1,\dots,P_N\}$, whose divisor is $D_1-D_2+(k-l)(\infty_+-\infty_-)$. 
$\qed$
\end{cor}

\subsection{The tau functions for the finite Toda lattice}\label{sec7}

It is classically known that the Toda equation admits a family of algebraic solutions which is expressed as a ratio of determinants.  
We demonstrate to derive them in terms of theta functions associated with $X$.

Let $\vect{v}_\pm(t)$ be the vector-valued polynomial function in $\lambda$ defined in \S \ref{sec2.1.1}.
We assign these vectors $\vect{v}_\pm(t)$ with a vector-valued function over $X$ in such a way that
(i) by replacing $\lambda$'s in each component of $\vect{v}_-(t)$ with $x$ and 
(ii) by replacing $\lambda$'s in each component of $\vect{v}_+(t)$ with $y$.
In this view point, the equation $F(t)\cdot \vect{v}_-(t)=\vect{v}_+(t)$ (Eq.\,(\ref{eq7})) represents the perversity of two vector-valued functions $\vect{v}_-(t)$, $\vect{v}_+(t)$ at the intersection $X_-\cap X_+$.

We regard the pair $\vect{v}(t)=(\vect{v}_-(t),\vect{v}_+(t))$ as one vector-valued function over $X\setminus\{P_1,\dots,P_N\}$.
Denote the $k$-th component of $\vect{v}(t)$ by $v_k$.
From (\ref{eq7}) and (\ref{eq8}), we have
$$
\phi_{k-1,N-k}(D_k(t))=\overline{F}(t)=[\exp(t\lambda_1):\cdots:\exp(t\lambda_N)]\in \mathrm{Jac}(X)\simeq (\CC^\times)^N/\CC^\times,
$$
where $D_k(t)=(v_k)_0\in S_{k-1,N-k}$.

\begin{prop}\label{prop10}
Let $1\leq k,l\leq N$.
Suppose that there exists two general divisors $D_k(t)\in S_{k-1,N-k}$, $D_l(t)\in S_{l-1,N-l}$ which satisfy the equation
$$
\phi_{k-1,N-k}(D_k(t))
=
\phi_{l-1,N-l}(D_l(t))
=
\overline{F}(t).
$$
Then, there exists a constant $C_{k,l}(t)$ such that
$$
C_{k,l}(t)\cdot \dfrac{v_k(p)}{v_l(p)}=
\dfrac{\Theta_k(p;F(t))}{\Theta_l(p;F(t))},\qquad p\in X\setminus\{P_1,\dots,P_N\}.
$$
\end{prop}
\proof 
Since $(v_k(p))=D_k(t)-(k-1)\infty_--(N-k)\infty_+$, the function $\frac{v_k(p)}{v_l(p)}$ is a single-valued function over $X$ whose divisor is $D_k(t)-D_l(t)+(k-l)(\infty_+-\infty_-)$.
On the other hand, by Corollary \ref{cor2.7}, the function $\frac{\Theta_k(p;F(t))}{\Theta_l(p;F(t))}$ possesses the same properties as well.
Therefore, they must coincides as rational functions over $X$ up to constant. $\qed$

\begin{lemma}\label{lemma11}
For $p\in X$, let $x=x(p)$ be the rational coordinate of $X_-$ and $y=y(p)$ be that of $X_+$.
The behaviors of $\dfrac{\Theta_k(p;F(t))}{\Theta_l(p;F(t))}$ when $p\to \infty_\pm$ are given by the followings:
\begin{gather*}
\dfrac{\Theta_k(p;F(t))}{\Theta_l(p;F(t))}\stackrel{p\to \infty_-}{\sim}
\dfrac{\tau_{k-1} x^{k-1}-\tau_{k-1}'x^{k-2}+\cdots}{\tau_{l-1} x^{l-1}-\tau_{l-1}'x^{l-2}+\cdots},\qquad
\dfrac{\Theta_k(p;F(t))}{\Theta_l(p;F(t))}\stackrel{p\to \infty_+}{\sim} 
\dfrac{\tau_{k} y^{N-k}+\cdots}{\tau_{l} y^{N-l}+\cdots},
\end{gather*}
where
$\vect{\lambda}^n:=(\lambda_1^n,\dots,\lambda_N^{n})^T$,
$\vect{F}:=(\exp(t\lambda_1),\dots,\exp(t\lambda_N))^T$ and
\begin{gather}
\tau_n:=\tau_n(t)=(-1)^{n-1}\cdot
\det (
\vect{1},\vect{\lambda},\dots,\vect{\lambda}^{N-n+1},
\vect{F},\vect{F}\vect{\lambda},\cdots,\vect{F}\vect{\lambda}^{n-1}
),\label{e18}\\
\tau_n':=\tau_n'(t)=(-1)^{n-1}\cdot
\det (
\vect{1},\vect{\lambda},\dots,\vect{\lambda}^{N-n+1},
\vect{F},\vect{F}\vect{\lambda},\cdots,\vect{F}\vect{\lambda}^{n-2},
\vect{F}\vect{\lambda}^{n}
).\label{e19}
\end{gather}
\end{lemma}
\proof The formulas are proved directly from (\ref{theta}).\qed
\begin{cor}
We have $C_{k,l}(t)=\dfrac{\tau_k}{\tau_l}$, hence we conclude
$\dfrac{\tau_k}{\tau_l}\cdot\dfrac{v_k(p)}{v_l(p)}=
\dfrac{\Theta_k(p;F(t))}{\Theta_l(p;F(t))}$,
$(\forall p\in X)$.
\end{cor}
\proof 
From $v_k(p)\sim x^{k-1}$ ($p\to \infty_-$), Proposition \ref{prop10} and Corollary \ref{lemma11}.
$\qed$

Let $l<k$.
When $p\to \infty_+$, it follows that $\dfrac{v_{k}(p)}{v_l(p)}\sim b_lb_{l+1}\cdots b_{k-1}y^{k-l}$.
Therefore, from Proposition \ref{prop10}, we obtain the equation
$$
b_lb_{l+1}\cdots b_{k-1}=
C_{k,l}(t)^{-1}\dfrac{\tau_{k}}{\tau_{l}}=\dfrac{\tau_{l-1}\tau_{k}}
{\tau_{k-1}\tau_{l}}.
$$
By substituting $l=n$, $k=n+1$, we have 
\begin{equation}\label{todab}
b_n=\dfrac{\tau_{n-1}\tau_{n+1}}{\tau_{n}^2}.
\end{equation}

$a_i$ is also given by the theta functions.
From (\ref{eq4}), we have $b_{n-1}v_{n-1}(p)+(a_{n}-\lambda)v_{n}(p)+v_{n+1}(p)=0$.
Then, $a_{n}=\lambda-b_{n-1}\dfrac{v_{n-1}(p)}{v_{n}(p)}-\dfrac{v_{n+1}(p)}{v_{n}(p)}$.
By taking the limit $p\to \infty_-$ (see Lemma \ref{lemma11}), we obtain
\begin{align*}
x-b_{n-1}\dfrac{\tau_{n-1}}{\tau_{n-2}}\dfrac{\Theta_{n-1}(p)}{\Theta_{n}(p)}-
\dfrac{\tau_{n-1}}{\tau_{n}}\dfrac{\Theta_{n+1}(p)}{\Theta_{n}(p)}
\sim x-
\dfrac
{
\tau_{n-1}}{\tau_{n}}\dfrac{\tau_{n}x^{n}
-\tau_{n}'x^{n-1}+\cdots
}
{
\tau_{n-1}x^{n-1}-\tau_{n-1}'x^{n-2}+\cdots
}
\mathop{\longrightarrow}^{x\to \infty}
\frac{\tau_{n}'}{\tau_{n}}-\frac{\tau_{n-1}'}{\tau_{n-1}},
\end{align*}
which implies
\begin{equation}\label{todaa}
a_{n}=\dfrac{\tau_{n}'}{\tau_{n}}-\dfrac{\tau_{n-1}'}{\tau_{n-1}}.
\end{equation}

\begin{thm}\label{kai1}
If the divisors $D_n$ are general for $1\leq n< N$, the variables $a_n=a_n(t)$, $b_n=b_n(t)$ in the Toda equation are given by (\ref{todab}) and (\ref{todaa}).
$\qed$
\end{thm}
In terms of the rational map 
\begin{equation}\label{kric}
\Phi:\mathcal{T}_\Lambda\to \mathrm{Jac}(X);\qquad L\mapsto \overline{F}(t)=F(t)\,\mathrm{mod}\,{\CC^\times},
\end{equation}
which is defined from (\ref{eq7}), Theorem \ref{kai1} is expressed as follows:
\begin{cor}
There exists a Zariski open set $\mathcal{U}\subset \mathrm{Jac}(X)$ and a rational map $\Psi:\mathcal{U}\to \mathcal{T}_\Lambda$; $\overline{F}(t)\mapsto \{a_n(t),b_n(t)\}$ such that $\Psi\circ\Phi=\mathrm{id}_{\mathcal{T}_\Lambda}$ and $\Phi\circ\Psi=\mathrm{id}_\mathcal{U}$.
Especially, $\Phi$ is injective. \qed
\end{cor}

%
%

\section{Totally non-negative part of the phase space}\label{sec3a}

In this section, we introduce some geometric characterization of the totally non-negative part (TNN part) of the phase space.

Denote $\mathcal{T}_\Lambda=H^{-1}(\gamma(\Lambda))$, the isospectral set associated with the spectrum set $\Lambda$ (see \S \ref{sec1}).
Let $X$ be a singular curve defined in the previous section.
In this section, we study the TNN part $\mathcal{T}_\Lambda^\geq\subset \mathcal{T}_\Lambda$ of the isospectral set.
One of the main results of this paper is the characterization of $\Phi(\mathcal{T}_\Lambda^\geq)$ as a subset of $\mathrm{Jac}(X)$.

\subsection{Properties of totally non-negative tridiagonal matrices}

We briefly introduce some of fundamental facts about TNN matrices.
For readers who are interested in the topic of TNN matrices, we recommend the text \cite{pinkus2010totally}.
We follow the notation of this text here.
\begin{defi}
A TNN matrix $L$ is said to be {\it irreducible} if $L^k$ is TP for some natural number $k$.
\end{defi}

\begin{thm}[Gantmacher-Krein]\label{gmt}
All the eigenvalues of an irreducible TNN matrix are simple and positive.
\end{thm}
\proof See \cite[Section 5]{pinkus2010totally}, for example. $\qed$

\begin{defi}
For an $N\times N$ matrix $L$ and series of integers
$1\leq i_1<i_2<\cdots<i_k\leq N$, $1\leq j_1<j_2<\cdots<j_k\leq N$,
we denote by
$
L
\left[ 
\begin{array}{c@{\,}c@{\,}c@{\,}c}
i_1 & i_2 & \cdots & i_k \\
j_1 & j_2 & \cdots & j_k 
\end{array}
\right] 
$
the $k\times k$ submarix of $L$ which consisted of $i_1,i_2,\dots,i_k$-th rows and $j_1,j_2,\dots,j_k$-th columns.
\end{defi}

\begin{prop}(\cite[Theorem 4.3]{pinkus2010totally}).\label{prop3.3}
An $N\times N$ tridiagonal matrix $L$ is TNN if and only if both of the following (i) and (ii) are satisfied:
(i)
For any $1\leq n\leq m\leq N$,
$
\det L
\left[ 
\begin{array}{c@{\ }c@{\,}c@{\,}c}
n & n+1 & \cdots & m \\
n & n+1 & \cdots & m 
\end{array}
\right] \geq 0
$.
(ii) Each off-diagonal entry ({\it i.e.}, an entry which is not a diagonal entry) of $L$ is nonnegative. $\qed$
\end{prop}

\begin{cor}\label{kei}
An $N\times N$ tridiagonal matrix $L$ is TNN if and only if all of the following (i)--(iii) are satisfied:
(i) $\det L\geq 0$, 
(ii) 
$L
\left[ 
\begin{array}{c@{\ }c@{\,}c@{\,}c}
2 & 3 & \cdots & N \\
2 & 3 & \cdots & N 
\end{array}
\right] $ is TNN,
(iii) Each off-diagonal component of $L$ is nonnegative.
\end{cor}
\proof 
Necessity is obvious.
We prove the sufficiency.
Suppose that $L$ satisfies (i)--(iii).
From Proposition \ref{prop3.3}, it suffices to prove
$\det L
\left[ 
\begin{array}{c@{\ }c@{\,}c@{\,}c}
1 & 2 & \cdots & k \\
1 & 2 & \cdots & k 
\end{array}
\right]
$ is TNN for any $k=1,2,\dots,N$.
For $k=N$, this is obvious from (i).
Assume $1\leq k\leq N-1$.
By Sylvester's relation, we have
\begin{align*}
&
\det L
\left[ 
\begin{array}{c@{\ }c@{\,}c@{\,}c}
1 & 2 & \cdots & k \\
1 & 2 & \cdots & k 
\end{array}
\right]
\cdot 
\det L
\left[ 
\begin{array}{c@{\ }c@{\,}c@{\,}c}
2 & 3 & \cdots & k+1 \\
2 & 3 & \cdots & k+1
\end{array}
\right]\\
&=
\det L
\left[ 
\begin{array}{c@{\ }c@{\,}c@{\,}c}
1 & 2 & \cdots & k+1 \\
1 & 2 & \cdots & k+1
\end{array}
\right]
\cdot 
\det L
\left[ 
\begin{array}{c@{\ }c@{\,}c@{\,}c}
2 & 3 & \cdots & k-1 \\
2 & 3 & \cdots & k-1
\end{array}
\right]+
\det L
\left[ 
\begin{array}{c@{\ }c@{\,}c@{\,}c}
1 & 2 & \cdots & k \\
2 & 3 & \cdots & k+1 
\end{array}
\right]
\cdot 
\det L
\left[ 
\begin{array}{c@{\ }c@{\,}c@{\,}c}
2 & 3 & \cdots & k+1 \\
1 & 2 & \cdots & k
\end{array}
\right].
\end{align*}
From (ii), it follows that
$
\det L
\left[ 
\begin{array}{c@{\ }c@{\,}c@{\,}c}
2 & 3 & \cdots & k-1 \\
2 & 3 & \cdots & k-1
\end{array}
\right]\geq 0$,
$\det L
\left[ 
\begin{array}{c@{\ }c@{\,}c@{\,}c}
2 & 3 & \cdots & k+1 \\
2 & 3 & \cdots & k+1
\end{array}
\right]\geq 0
$.
Since $L$ is tridiagonal and (iii), we have
$
\det L
\left[ 
\begin{array}{c@{\ }c@{\,}c@{\,}c}
1 & 2 & \cdots & k \\
2 & 3 & \cdots & k+1
\end{array}
\right]\geq 0$,
$\det L
\left[ 
\begin{array}{c@{\ }c@{\,}c@{\,}c}
2 & 3 & \cdots & k+1 \\
1 & 2 & \cdots & k
\end{array}
\right]\geq 0
$.
Finally, we conclude that
$
\det L
\left[ 
\begin{array}{c@{\ }c@{\,}c@{\,}c}
1 & 2 & \cdots & k+1 \\
1 & 2 & \cdots & k+1
\end{array}
\right]\geq 0
$ implies
$
\det L
\left[ 
\begin{array}{c@{\ }c@{\,}c@{\,}c}
1 & 2 & \cdots & k \\
1 & 2 & \cdots & k
\end{array}
\right]\geq 0
$.
Therefore, by (i), we conclude that
$
\det L
\left[ 
\begin{array}{c@{\ }c@{\,}c@{\,}c}
1 & 2 & \cdots & k \\
1 & 2 & \cdots & k
\end{array}
\right]
$ is TNN for any $k$.
$\qed$
\begin{rem}
Corollary \ref{kei} is also true if we replace the condition (ii) with the following (ii)$'$:
\begin{center}
(ii)$'$
$L
\left[ 
\begin{array}{c@{\ }c@{\,}c@{\,}c}
1 & 2 & \cdots & N-1 \\
1 & 2 & \cdots & N-1 
\end{array}
\right] $ is TNN.
\end{center}
\end{rem}

\begin{thm}\label{thm3.6}
Let $N\geq 2$.
For $N\times N$ irreducible TNN matrix $L$, set
$Q:=L\left[ 
\begin{array}{c@{\ }c@{\,}c@{\,}c}
2 & 3 & \cdots & N \\
2 & 3 & \cdots & N 
\end{array}
\right]$,
$Q':=L\left[ 
\begin{array}{c@{\ }c@{\,}c@{\,}c}
1 & 2 & \cdots & N-1 \\
1 & 2 & \cdots & N-1 
\end{array}
\right]$.
Let $0<\lambda_1<\lambda_2<\dots<\lambda_N$ be the spectrum of $L$,
$0<\mu_1<\mu_2<\dots<\mu_{N-1}$ be the spectrum of $Q$ and $0<\mu'_1<\mu'_2<\dots<\mu'_{N-1}$ be the spectrum of $Q'$.
Then,
\begin{gather}
0<\lambda_1<\mu_1<\lambda_2<\mu_2<\lambda_3<\cdots<\mu_{N-1}<\lambda_N,\label{e26}
\\
0<\lambda_1<\mu'_1<\lambda_2<\mu'_2<\lambda_3<\cdots<\mu'_{N-1}<\lambda_N
\label{e26a}\tag{$\ref{e26}'$}.
\end{gather}
\end{thm}
\proof \cite[Proposition 5.4]{pinkus2010totally}. $\qed$

\begin{prop}\label{prop3.7}
Assume $N\geq 2$.
Suppose that $L$ is a tridiagonal matrix whose off-diagonal entries are positive.
Let $Q$, $Q'$ be the same matrices as in Theorem \ref{thm3.6},
$\lambda_1<\dots<\lambda_N$ the spectrum of $L$, $\mu_1,\dots,\mu_{N-1}$ the spectrum of $Q$ and
$\mu'_1,\dots,\mu'_{N-1}$ the spectrum of $Q'$.
Then, the followings are equivalent:
(i) $L$ is TNN,
(ii) (\ref{e26}) holds,
(iii) (\ref{e26a}) holds.
\end{prop}
\proof We prove by induction on $N\geq 2$.
The case of $N=2$ is directly proven.
Assume $N\geq 3$ and $\lambda_i,\mu_i,\mu_i'$ to satisfy (\ref{e26}).
Let $L_\lambda:=L-\lambda E$,
$
X(\lambda):=\det L_\lambda$,
$X_1(\lambda):=
\det
L_\lambda\left[ 
\begin{array}{c@{\ }c@{\,}c@{\,}c}
2 & 3 & \cdots & N \\
2 & 3 & \cdots & N 
\end{array}
\right] $,
$
X_N(\lambda):=
\det
L_\lambda\left[ 
\begin{array}{c@{\ }c@{\,}c@{\,}c}
1 & 2 & \cdots & N-1 \\
1 & 2 & \cdots & N-1 
\end{array}
\right]
$,
$X_{1,N}(\lambda):=
\det
L_\lambda\left[ 
\begin{array}{c@{\ }c@{\,}c@{\,}c}
2 & 3 & \cdots & N-1 \\
2 & 3 & \cdots & N-1 
\end{array}
\right] $,
$Y(\lambda):=
\det
L_\lambda\left[ 
\begin{array}{c@{\ }c@{\,}c@{\,}c}
2 & 3 & \cdots & N \\
1 & 2 & \cdots & N-1 
\end{array}
\right] $ and
$
Z(\lambda):=
\det
L_\lambda\left[ 
\begin{array}{c@{\ }c@{\,}c@{\,}c}
1 & 2 & \cdots & N-1 \\
2 & 3 & \cdots & N 
\end{array}
\right]
$.
From Silvester's relation (see \cite[(p.5, Eq.\,(1.2))]{pinkus2010totally}, for example), we have the equation
\begin{equation}
X(\lambda)X_{1,N}(\lambda)=X_1(\lambda)X_N(\lambda)-Y(\lambda)Z(\lambda).
\end{equation}
Since $L$ is tridiagonal and its off-diagonal components are positive, $Y(\lambda)$ and
$Z(\lambda)$ are positive constants.
Then,
$
X(\lambda)X_{1,N}(\lambda)<X_1(\lambda)X_N(\lambda)
$.
Substituting
$\lambda=\mu_i$ ($i=1,\dots,N-1$), we obtain $X(\mu_i)X_{1,X}(\mu_i)<0$.
From (\ref{e26}), the inequalities 
$$
X(0)=\det L>0,\quad 
X(\mu_1)<0,\quad 
X(\mu_2)>0,\quad 
X(\mu_3)<0,\dots
$$ 
hold, which imply 
$$
X_{1,X}(\mu_1)>0,\quad
X_{1,X}(\mu_2)<0,\quad
X_{1,X}(\mu_3)>0,\dots.
$$
By the intermediate value theorem, there exists some real number $\mu_i<c_i<\mu_{i+1}$, $(i=1,\dots,N-2)$ with
$X_{1,N}(c_i)=0$.
As $X_{1,N}(\lambda)$ is a polynomial in $\lambda$ of degree $(N-2)$, all the roots of $X_{1,N}(\lambda)$ are expressed as $\lambda=c_1,\dots,c_{N-2}$.
By seeing the $(N-2)\times (N-2)$ matrix
$L
\left[ 
\begin{array}{c@{\ }c@{\,}c@{\,}c}
2 & 3 & \cdots & N-1 \\
2 & 3 & \cdots & N-1 
\end{array}
\right]
$ as a submatrix of $Q$, we conclude that $Q$ is TNN by hypothesis of induction.
From Corollary \ref{kei}, also $L$ is TNN. $\qed$

\subsection{Characterization of the TNN part}\label{sec3.2}

From Theorem \ref{gmt}, $\mathcal{T}_\Lambda^\geq \neq \emptyset$ implies that $\Lambda$ consists of distinct $N$ positive numbers.
Hereafter, we assume $\Lambda=\{\lambda_1,\dots,\lambda_N\}$ ($0<\lambda_1<\dots<\lambda_N$).
Let $X$ be the singular curve defined in the previous section and 
$
\mathrm{Jac}(X)=(\CC^\times)^N/\CC^\times 
$
be the generalized Jacobi variety of $X$.

The following theorem gives the characterization of the TNN part of the Toda flow.
\begin{thm}\label{thm3.8}
The image $\Phi(\mathcal{T}_\Lambda^\geq)\subset \mathrm{Jac}(X)$ of the TNN part is expressed as
\[
\Phi(\mathcal{T}_\Lambda^\geq)=
\left\{
[F_1:F_2:\dots:F_N]\in\mathrm{Jac}(X)\,\left\vert\, (-1)^iF_i>0\mbox{ for all } i
\right.
\right\}.
\]
\end{thm}
\proof 
Set 
$
M=
\left\{
[F_1:F_2:\dots:F_N]\in\mathrm{Jac}(X)\,\left\vert\, (-1)^iF_i>0\mbox{ for all } i
\right.
\right\}
$.

(Proof of $\Phi(\mathcal{T}_\Lambda^\geq)\subset M$.)
Suppose $L\in\mathcal{T}^\geq_\Lambda$.
From the construction of $\Phi$, we have
$
\Phi(L)=
[
X_1(\lambda_1):X_1(\lambda_2):\cdots:X_1(\lambda_N)
]
$,
where
$
X_1(\lambda)=
\det 
L_\lambda
\left[ 
\begin{array}{c@{\ }c@{\,}c@{\,}c}
2 & 3 & \cdots & N \\
2 & 3 & \cdots & N 
\end{array}
\right]
$.
By Theorem \ref{thm3.6}, we have $(-1)^{i-1}X_1(\lambda_i)>0$ for all $i$ which implies $\Phi(L)\in M$.

(Proof of $\Phi(\mathcal{T}_\Lambda^\geq)\supset M$.)
Let $[F_1:\dots:F_N]\in M$. 
Without loss of generality, we assume $(-1)^{i-1}F_i>0$ for all $i$.
Define $\tau_n$ and $\tau_n'$ by equations (\ref{e18}), (\ref{e19}) and set
\[
b_n:=\frac{\tau_{n-1}\tau_{n+1}}{\tau_{n}^2},\qquad
a_n:=\frac{\tau_{n}'}{\tau_{n}}-\frac{\tau_{n-1}'}{\tau_{n-1}}.
\] 
Let $L$ be the Lax matrix (\S 2.1) defined by $\{a_n,b_n\}$.
Then it follows that $\Phi(L)=[F_1:\dots:F_N]$.
Let us prove $L\in \mathcal{T}_\Lambda^\geq$.
For a series of natural numbers $1\leq i_1<i_2<\cdots<i_{n}\leq N$, set
$
\Delta_{i_1,i_2,\dots,i_{n}}:=
\prod_{1\leq a<b\leq n}
(\lambda_{i_b}-\lambda_{i_a})
$. 
Denote 
$$
\Delta^\dagger_{i_1,i_2,\dots,i_{n}}:=\Delta_{j_1,j_2,\dots,j_{N-n}},
$$
where $\{1,2,\dots,N\}=\{i_1,i_2,\dots,i_{n}\}
\sqcup
\{j_1,j_2,\dots,j_{N-n}\}
$.
By the Laplace expansion, $\tau_n$ is expanded as follows:
\[
\tau_n=\sum_{1\leq i_1<i_2<\cdots<i_{n-1}\leq N}
(-1)^{i_1+i_2+\cdots+i_{n-1}}\cdot 
F_{i_1}F_{i_2}\cdots F_{i_{n-1}}\cdot
\Delta_{i_1,\dots,i_{n-1}}\cdot \Delta^\dagger_{i_1,\dots,i_{n-1}}.
\]
By the inequalities $0\leq \lambda_1<\cdots<\lambda_N$ and $(-1)^{i-1}F_i>0$, all the numbers
$\Delta_{i_1,\dots,i_{n-1}}$, $\Delta^\dagger_{i_1,\dots,i_{n-1}}$, $\tau_n$, $b_n$ must be positive.
Since $b_n$ is a ratio of $\tau_m$'s, the Lax matrix $L$ is a matrix whose off-diagonal components are positive.
Let $X_1(\lambda)$ be the $(1,1)$-cofactor of the matrix $L-\lambda E$.
By definition of the Abel-Jacobi map, there exists some complex number $c$ such that $X_1(\lambda_i)=c\cdot F_i$, ($i=1,\dots,N$).
Since $X_1(\lambda)$ is a real polynomial of degree $(N-1)$ whose leading coefficient is $(-1)^{N-1}$, $c$ must be real and positive.
Hence, $(-1)^{i-1}X_1(\lambda_i)>0$ ($\forall i$).
By the intermediate value theorem, there exists a real number $\lambda_i<\mu_i<\lambda_{i+1}$ with $X_1(\mu_i)=0$.
By applying Proposition \ref{prop3.7} to $L$, the tridiagonal matrix whose off-diagonal components are positive, we conclude that $L$ is TNN. 
Then, $L\in \mathcal{T}_\Lambda^\geq$.
$\qed$

\subsection{The real structure}\label{sec3.3}

Let $\mathcal{T}_\Lambda^\RR:=\mathcal{T}_\Lambda\cap M_N(\RR)$ be the real part of $\mathcal{T}$ and 
$$
\mathrm{Jac}(X)_\RR:=\{[F_1:\dots:F_N]\,\vert\,F_i\in \RR\setminus\{0\}  \}
$$ 
be that of $\mathrm{Jac}(X)$.
Since the construction of $\Phi$ is closed in the reals, the linearization map $\Phi:\mathcal{T}_\Lambda\to \mathrm{Jac}(X)$ induces the algebraic map 
\[
\Phi_\RR:\mathcal{T}_\Lambda^\RR\to\mathrm{Jac}(X)_\RR 
\]
of real algebraic varieties.

The real part $\mathrm{Jac}(X)_\RR$ consists of $2^{N-1}$ connected components of dimension $N-1$.
Let $\mathrm{Jac}(X)_\RR^0\subset \mathrm{Jac}(X)_\RR^0$ be the connected component which contains the identity element $[1:1:\cdots:1]$.
Theorem \ref{thm3.8} states the existence of the isomorphism 
$$
\Phi_\geq:\mathcal{T}_\Lambda^\geq\to\mathrm{Jac}(X)_\RR^0
$$ 
of two semi-algebraic varieties.

\section*{Acknowledgments} 

The two of the authors (S.\,I. and K.\,N.) were partially supported by KAKENHI (25610008).
S.\,I. was also supported by KAKENHI (26800062).

\appendix

\bibliography{test2}


\end{document}